\documentclass[12pt]{article}

\usepackage{amsmath,amssymb}

\begin{document}

\pagestyle{myheadings} \markright{PRIME GEODESIC THEOREM...}

\title{A prime geodesic theorem for higher rank spaces}
\author{Anton Deitmar}

\date{}
\maketitle

\def \1{{\bf 1}}
\def \a{{{\mathfrak a}}}
\def \ad{{\rm ad}}
\def \al{\alpha}
\def \ar{{\alpha_r}}
\def \A{{\mathbb A}}
\def \Ad{{\rm Ad}}
\def \Aut{{\rm Aut}}
\def \b{{{\mathfrak b}}}
\def \bs{\backslash}
\def \B{{\cal B}}
\def \c{{\mathfrak c}}
\def \cent{{\rm cent}}
\def \C{{\mathbb C}}
\def \CA{{\cal A}}
\def \CB{{\cal B}}
\def \CC{{\cal C}}
\def \CD{{\cal D}}
\def \CE{{\cal E}}
\def \CF{{\cal F}}
\def \CG{{\cal G}}
\def \CH{{\cal H}}
\def \CHC{{\cal HC}}
\def \CL{{\cal L}}
\def \CM{{\cal M}}
\def \CN{{\cal N}}
\def \CP{{\cal P}}
\def \CQ{{\cal Q}}
\def \CO{{\cal O}}
\def \CS{{\cal S}}
\def \CT{{\cal T}}
\def \CV{{\cal V}}
\def \CW{{\cal W}}
\def \d{{\mathfrak d}}
\def \df{\ \begin{array}{c} _{\rm def}\\ ^{\displaystyle =}\end{array}\ }
\def \det{{\rm det}}
\def \diag{{\rm diag}}
\def \dist{{\rm dist}}
\def \End{{\rm End}}
\def \eqn{\begin{eqnarray*}}
\def \endeqn{\end{eqnarray*}}
\def \eps{\varepsilon}
\def \F{{\mathbb F}}
\def \Fx{{\mathfrak x}}
\def \FX{{\mathfrak X}}
\def \g{{{\mathfrak g}}}
\def \ga{\gamma}
\def \Ga{\Gamma}
\def \Gal{{\rm Gal}}
\def \h{{{\mathfrak h}}}
\def \Hom{{\rm Hom}}
\def \im{{\rm im}}
\def \Im{{\rm Im}}
\def \ind{{\rm ind}}
\def \k{{{\mathfrak k}}}
\def \K{{\cal K}}
\def \l{{\mathfrak l}}
\def \la{\lambda}
\def \lap{\triangle}
\def \li{{\rm li}}
\def \La{\Lambda}
\def \Lie{{\rm Lie}}
\def \m{{{\mathfrak m}}}
\def \mod{{\rm mod}}
\def \n{{{\mathfrak n}}}
\def \name{\bf}
\def \Mat{{\rm Mat}}
\def \N{\mathbb N}
\def \o{{\mathfrak o}}
\def \ord{{\rm ord}}
\def \O{{\cal O}}
\def \p{{{\mathfrak p}}}
\def \ph{\varphi}
\def \prf{\noindent{\bf Proof: }}
\def \Per{{\rm Per}}
\def \q{{\mathfrak q}}
\def \qed{\ifmmode\eqno $\square$\else\noproof\vskip 12pt plus 3pt minus 9pt \fi}
 \def\noproof{{\unskip\nobreak\hfill\penalty50\hskip2em\hbox{}%
     \nobreak\hfill $\square$\parfillskip=0pt%
     \finalhyphendemerits=0\par}}
\def \Q{\mathbb Q}
\def \res{{\rm res}}
\def \R{{\mathbb R}}
\def \Re{{\rm Re \hspace{1pt}}}
\def \r{{\mathfrak r}}
\def \ra{\rightarrow}
\def \rank{{\rm rank}}
\def \SL{{\rm SL}}
\def \SO{{\rm SO}}
\def \supp{{\rm supp}}
\def \Spin{{\rm Spin}}
\def \t{{{\mathfrak t}}}
\def \T{{\mathbb T}}
\def \tr{{\hspace{1pt}\rm tr\hspace{2pt}}}
\def \vol{{\rm vol}}
\def \z{\zeta}
\def \Z{\mathbb Z}
\def \={\ =\ }

\newcommand{\frack}[2]{\genfrac{}{}{0pt}{}{#1}{#2}}
\newcommand{\rez}[1]{\frac{1}{#1}}
\newcommand{\der}[1]{\frac{\partial}{\partial #1}}
\renewcommand{\binom}[2]{\left( \begin{array}{c}#1\\#2\end{array}\right)}
\newcommand{\norm}[1]{\parallel #1 \parallel}
\renewcommand{\matrix}[4]{\left(\begin{array}{cc}#1 & #2 \\ #3 & #4 \end{array}\right)}
\renewcommand{\sp}[2]{\langle #1,#2\rangle}
\renewcommand{\labelenumi}{(\alph{enumi})}

\newtheorem{theorem}{Theorem}[section]
\newtheorem{conjecture}[theorem]{Conjecture}
\newtheorem{lemma}[theorem]{Lemma}
\newtheorem{corollary}[theorem]{Corollary}
\newtheorem{proposition}[theorem]{Proposition}

{\bf Abstract.}
{\small A prime geodesic theorem for regular geodesics in a higher rank locally
symmetric space is proved. An application
to class numbers is given.  
The proof relies on a
Lefschetz formula for higher rank torus actions. }  
$ $

\tableofcontents

\newpage

\begin{center} {\bf Introduction} \end{center}

The prime geodesic
theorem gives a growth asymptotic for the number of closed geodesics counted by their lengths
\cite{Gangolli,Hejhal,Knieper,Koyama,LuoSarnak,Margulis,Pollicott,Zelditch}.
It has hitherto only
been proven for manifolds of strictly negative curvature. 
For manifolds containing higher dimensional flats it is not a priori clear what a
prime geodesic theorem might look like. 
In this paper we give such a
theorem for locally symmetric spaces of arbitrary rank, i.e., they may contain
higher dimensional flats. The regular geodesics in such a space give points in a
higher dimensional Weyl cone. Therefore it is natural to expect a prime geodesic
theorem which gives asymptotics in several variables.
 We
introduce a new  analytical function which could be viewed as higher
dimensional analogue of the logarithmic derivative of
the Selberg zeta function. Of this we 
get analytic continuation only in a certain translate of the
positive Weyl chamber, but this actually suffices to prove
the prime geodesic theorem.

We describe the main result of the paper. One of the
various equivalent formulations of the prime geodesic theorem for locally symmetric spaces of rank one is
the following. Let $\bar X$ be a compact locally symmetric
space with universal covering of rank one. For $T>0$ let
$$
\psi(T)\=\sum_{c\,:\, e^{l(c)}\le T} l(c_0).
$$
Here the sum runs over all closed geodesics $c$ such that
$e^{l(c)}\le T$, where $l(c)$ is the length of the
geodesic $c$, and $c_0$ is the prime geodesic underlying
$c$. Then, under a suitable scaling of the metric, as $T\ra\infty$,
$$
\psi(T)\ \sim\ T.
$$

We now replace the space $\bar X$ by an arbitrary compact locally symmetric
space which is a quotient of a globally symmetric space $X=G/K$ where $G$ is a semisimple
Lie group of split-rank
$r$ and
$K$ a maximal compact subgroup.  A closed geodesic $c$ gives rise to
a point $a_c$ in the closure of the negative Weyl chamber $A^-$ of a maximal split torus
$A$.  We assume that the geodesic is regular \cite{jost}, i.e., the element $a_c$ lies in the
interior of the Weyl chamber. (In rank one every geodesic is regular.) For
$T_1,\dots T_r>0$ let
$$
\psi(T_1,\dots T_r)\=\sum_{c\,:\, a_{c,j}\le T_j} \la_c,
$$
where $\la_c$ is the volume of the unique maximal flat $c$ lies in and $a_{c,j}$ are the coordinates of $a_c$ with
respect to a canonical coordinate system given by the roots. The sum runs over all
regular closed geodesics modulo homotopy. The main result of this paper is that,
as $T_j$  tends to infinity for every $j$,
$$
\psi(T_1,\dots T_r)\ \sim\ T_1\cdots T_r.
$$

The  proof is based on a Lefschetz formula that generalizes the pioneering work of Andreas
Juhl
\cite{juhl}. The formula can be interpreted as a dynamical Lefschetz formula as follows. The
choice of a negative Weyl chamber $A^-$ induces a stable foliation $\CF$. The torus $A$
acts on the reduced tangential cohomology $\bar H^\bullet(\CF)$ of this foliation and the
Lefschetz formula states that there is an identity of distributions on the negative Weyl
chamber
$A^-$,
$$
\sum_{q=0}^{\rank\, \CF} (-1)^q\, \tr(a\mid \bar H^q(\CF))\=\sum_c\ind(c)\delta(a_c),
$$
where the sum on the right hand side runs over all regular closed geodesics,
$\delta(a_c)$ is the $\delta$-distribution at $a_c$ and $\ind(c)$ is a certain index which
incorporates monodromy data of $c$.

A specific test function is constructed so that the local side of the Lefschetz formula gives a
high derivative of the generalized Dirichlet series $L(s)=\sum_c \ind(c) a_c^s$. The Lefschetz
formula allows for a detailed analysis of this derivative of $L(s)$. By methods of analytical
number theory one finally gets the prime geodesic theorem. 

The methods of this paper will not give a similar result for non-regular geodesics. Though there is a Lefschetz formula
for non-regular geodesics, in general it contains Euler characteristics which may be positive or negative. The
Wiener-Ikehara Theorem used, however, relies on positivity. 

In an appendix we give an application to class numbers in totally real number fields.

\section{The Lefschetz formula}\label{sec1}
In this section we prove a Lefschetz formula for higher rank groups.  We give a more
general version then actually needed in section \ref{sec2}, as we allow twists by
nontrivial $M$-representations $\sigma$ (see below).

Let $G$ be a connected semisimple Lie group with finite
center and choose a maximal compact subgroup $K$ with
Cartan involution $\theta$, i.e., $K$ is the group of
fixed points of $\theta$. Let $P$ be a minimal parabolic subgroup with Langlands decomposition
$P=MAN$. Modulo conjugation we can assume that $A$ and $M$ are stable under
$\theta$. Then $A$ is a maximal split torus of $G$ and $M$ is a subgroup of $K$. The centralizer of $A$
is $AM$. Let $W(A,G)$ be the \emph{Weyl group} of $A$, i.e. $W$ is the quotient of the
normalizer of $A$ by the centralizer. This is a finite group acting on $A$. 

We have to fix Haar measures. We use the normalization of
Harish-Chandra \cite{HC-HA1}. Note that this normalization
depends on the choice of an invariant bilinear form $B$ on
$\g_\R$ which we keep at our disposal until later. Changing $B$ amounts to scaling the
metric of the symmetric space. Note further that in this normalization of Haar measures the
compact groups $K$ and $M$ have volume $1$.

We write $\g_\R, \k_\R,\a_\R,\m_\R,\n_\R$ for the real Lie
algebras of $G,K,A,M,N$ and $\g,\k,\a,\m,\n$ for their
complexifications. $U(\g)$ is the universal enveloping
algebra of $\g$. This algebra is isomorphic to the algebra
of all left invariant differential operators on $G$ with
complex coefficients. Pick a Cartan subalgebra $\t$ of $\m$. Then $\h=\a\oplus\t$ is a
Cartan subalgebra of $\g$. Let $W(\h,\g)$ be the corresponding absolute Weyl group.

Let $\a^*$ denote the dual space of the complex vector
space $\a$. Let $\a_\R^*$ be the real dual of $\a_\R$. We identify $\a_\R^*$ with the real vector space of all $\la\in\a^*$ that map $\a_\R$ to $\R$. Let $\Phi\subset\a^*$ be the set of all roots
of the pair $(\a,\g)$ and let $\Phi^+$ be the subset of
positive roots with respect to $P$. Let $\Delta\subset
\Phi^+$ be the set of simple roots. Then $\Delta$ is a basis of $\a^*$. The open {\it
negative Weyl chamber} $\a_\R^-\subset\a_\R$ is the cone of all
$X\in\a_\R$ with $\alpha(X)<0$ for every $\alpha\in\Delta$. Let $\overline{\a_\R^-}$ be the
closure of $\a_\R^-$. The $W(A,G)$-translates $w\a_\R^-$ of $\a_\R^-$ are pairwise disjoint
and their union equals $\a_\R$ minus a finite number of hyperplanes. This is called the
\emph{regular set},
$$
\a_\R^{reg}\df \bigcup_{w\in W(A,G)} w\a_\R^-.
$$
Let $A^{reg}=\exp\left(\a_\R^{reg}\right)$ be the regular set in $A$. The elements are
called the regular elements of $A$. A given $a\in A$ lies in $A^{reg}$ if and only if the
centralizer of $a$ in $G$ equals the centralizer of $A$ in $G$ which is $AM$.

The bilinear form $B$ is indefinite on $\g_\R$, but the form
$$
\sp{X}{Y}\df -B(X,\theta(Y))
$$
is positive definite, ie an inner product on $\g_\R$. We extend it to an inner product on
the complexification $\g$. Let $\norm X=\sqrt{\sp XX}$ be the corresponding norm.
The form $B$, being nondegenerate, identifies $\g$ to its dual space $\g^*$. In this way we
also define an inner product $\sp ..$ and the corresponding norm on $\g^*$.
Furthermore, if $V\subset\g$ is any subspace on which $B$ is nondegenerate, then $B$ gives
an identification of $V^*$ with $V$ and so one gets an inner product and a norm on
$V^*$. This in particular applies to $V=\h$, a Cartan subalgebra of $\g$, which is defines
over
$\R$.

Let $\Ga\subset G$ be a discrete, cocompact, torsion-free subgroup. 
We are interested in the closed geodesics on the locally symmetric space $X_\Ga=\Ga\bs
X=\Ga\bs G/K$. Every such geodesic 
$c$ lifts to a $\Ga$-orbit of geodesics on $X$ and gives a $\Ga$-conjugacy class $[\ga_c]$ of
elements closing the particular geodesics. This induces a bijection between the set of all
homotopy classes of closed geodesics in $X_\Ga$ and the set of all non-trivial conjugacy
classes in $\Ga$.

An element $am$ of $AM$ is called {\it split regular} if
$a$ is regular in $A$. An element $\ga$ of $\Ga$
is called split regular if $\ga$ is in $G$ conjugate to a
split regular element $a_\ga m_\ga$ of $AM$. In that case
we may (and do) assume that $a_\ga$ lies in the negative
Weyl chamber $A^-=\exp(\a_\R^-)$ in $A$. Let $\CE(\Ga)$
denote the set of all conjugacy classes in $\Ga$ that
consist of split regular elements. Via the above correspondence the set $\CE(\Ga)$ can be
identified with the set of all homotopy classes of regular closed geodesics in $X_\Ga$, \cite{jost}.

Let $[\ga]\in\CE(\Ga)$.
There is a closed geodesic $c$ in the Riemannian manifold $\Ga\bs G/K$ which gets closed by
$\ga$. This means that there is a lift $\tilde c$ to the universal covering $G/K$ which is
preserved by $\ga$ and $\ga$ acts on $\tilde c$ by a translation.
The closed geodesic $c$ is not unique in general but lies in a unique maximal flat
submanifold $F_\ga$ of $\Ga\bs G/K$. Let $\la_\ga$ be the volume of that flat,
$$
\la_\ga\df \vol(F_\ga).
$$ 
Let $\Ga_\ga$ and $G_\ga$ denote the centralizers of $\ga$ in $\Ga$ and $G$ respectively.
The flat $F_\ga$ is the image of $\Ga_\ga\bs G_\ga /K_\ga$ in $\Ga\bs G/K$, where we assume
that $K$ is chosen so that $K_\ga=G_\ga\cap K$ is a maximal compact subgroup of $G_\ga$. 
In Harish-Chandra's normalization $G_\ga$ is equipped with the Haar measure that satisfies
$\int_{G_\ga}=\int_{G_\ga/K_\ga}\int_{K_\ga}$, where $K_\ga$ has the Haar measure with
$\vol(K_\ga)=1$ and $G_\ga/K_\ga$ gets the measure induced by the metric of $G/K$. Therefore
$$
\la_\ga\= \vol(\Ga_\ga\bs G_\ga /K_\ga)\=\vol(\Ga_\ga\bs G_\ga).
$$
Note that for not maximally split elements $\ga$ the relation between $\la_\ga$ and
$\vol(\Ga_\ga\bs G_\ga)$ is more complicated \cite{geom,hr}.

Let $\n$ denote the complexified Lie algebra of $N$. For
any $\n$-module $V$ let $H_q(\n,V)$ and $H^q(\n,V)$ for $q=0,\dots,\dim\n$
be the Lie algebra homology and cohomology \cite{BorWall}.  Let $\hat G$ denote the unitary dual of $G$, i.e., the set of isomorphism classes of irreducible unitary representations of $G$.  For $\pi\in\hat G$ let $\pi_K$ be the $(\g,K)$-module of $K$-finite vectors. If
$\pi\in\hat G$, then $H_q(\n,\pi_K)$ and $H^q(\n,\pi_K)$ are finite
dimensional $AM$-modules \cite{HeSch}. Note that they a priori only are $(\a\oplus\m,M)$-modules, but since $A$ is isomorphic to its Lie algebra they are  $AM$-modules.

Note that $AM$ acts on the Lie algebra $\n$ of $N$ by the
adjoint representation. Let $[\ga]\in\CE(\Ga)$. Since
$a_\ga\in A^-$ it follows that every eigenvalue of $a_\ga
m_\ga$ on $\n$ is of absolute value $<1$. Therefore
$\det(1-a_\ga m_\ga | \n)\ne 0$.

For $[\ga]\in\CE(\Ga)$ let
$$
\ind(\ga)\=\frac{\la_\ga }{\det(1-a_\ga m_\ga \mid \n)}\ >\ 0.
$$
Since $\Ga$ is cocompact, the unitary $G$-representation on $L^2(\Ga\bs G)$ splits
discretely with finite multiplicities
$$
L^2(\Ga\bs G)\= \bigoplus_{\pi\in\hat G} N_\Ga(\pi)\pi,
$$
where $N_\Ga(\pi)$ is a non-negative integer and $\hat G$ is the unitary dual of $G$. Fix a finite dimensional irreducible representation
$\sigma$ of $M$ and denote by $\breve\sigma$ the dual representation. A {\it quasi-character} of $A$ is a
continuous group homomorphism to $\C^\times$. Via
differentiation the set of quasi-characters can be
identified with the dual space $\a^*$. For $\la\in\a^*$ we write $a\mapsto a^\la$ for the
corresponding quasicharacter on $A$. We denote by
$\rho\in\a^*$ the modular shift with respect to $P$, i.e., for $a\in A$ we have
$\det(a|\n)=a^{2\rho}$.

For a complex vector space $V$ on which $A$ acts linearly and $\la\in\a^*$ let $(V)_\la$
denote the generalized $(\la+\rho)$-eigenspace, i.e.,
$$
(V)_\la\= \{ v\in V\mid (a-a^{\la+\rho}Id)^n v=0\ \ {\rm for\ some\ } n\in\N\}.
$$
Since $H^p(\n,\pi_K)$ is finite dimensional, the Jordan Normal Form Theorem implies that
$$
H^p(\n,\pi_K) \= \bigoplus_{\nu\in\a^*} H^p(\n,\pi_k)_\nu.
$$
Let $T$ be a Cartan subgroup of $M$ and let $\t$ be its complex Lie algebra. Then $AT$ is a
Cartan subgroup of $G$. Let $\Lambda_\pi\in(\a\oplus\t)^*$ be a representative of the
infinitesimal character of $\pi$. By Corollary 3.32 of \cite{HeSch} it follows,
$$
H_p(\n,\pi_K)\=\bigoplus_{\nu=w\Lambda_{\pi}|_\a}H_p(\n,\pi_K)_\nu,
$$
where $w$ ranges over $W(\g,\h)$.

\begin{lemma}\label{1.1}
For $0\le p\le d=\dim(\n)$ we have
$$
 H_p(\n,\pi_K)\ \cong\ H^{d-p}(\n,\pi_K)\otimes\det(\n),
$$
where the determinant of a
finite dimensional space is the top exterior power. So $\det(\n)$ is a one dimensional
$AM$-module on which $AM$ acts via the quasi-character $am\mapsto \det(am|\n)=a^{2\rho}$.
This in particular implies
$$
H^p(\n,\pi_K)\=\bigoplus_{\nu= w\Lambda_\pi |_\a}H^p(\n,\pi_K)_{\nu-2\rho}.
$$
\end{lemma}
\prf The first part follows straight from the definition of Lie algebra cohomology. The
second part by Corollary 3.32 of \cite{HeSch}.
\qed

For $\la\in\a^*$
 and $\pi\in\hat G$ let
$$
m_\la^\sigma(\pi)\= \sum_{q=0}^{\dim \n} (-1)^{q+\dim\n} \dim\left(
H^q(\n,\pi_K)_\la\otimes\breve\sigma\right)^M,
$$
where the superscript $M$ indicates the subspace of $M$-invariants. Then 
$m_\la^\sigma(\pi)$ is an integer and by the above, the set of $\la$ for which
$m_\la^\sigma(\pi)\ne 0$ for a given $\pi$ has at most $|W(\g,\h)|$ many elements.

For  $\mu\in\a^{*}$ and $j\in\N$ let $\CC^{j,\mu,-}(A)$ denote the space of functions
$\ph$ on $A$ which \nopagebreak
\begin{itemize}
  \item are $j$-times continuously differentiable on $A$,
  \item are zero outside $A^-$,
  \item are such that $a^{-\mu} D\ph(a)$ is bounded on $A$  for every invariant differential operator $D$ on $A$ of
degree $\le j$.
\end{itemize}

For every invariant differential operator $D$ of degree $\le j$ let $N_D(\ph)=\sup_{a\in A}\left| a^{-\mu}\, D\ph(a)\right|$.
Then $N_D$ is a seminorm.
Let $D_1,\dots,D_n$ be a basis of the space of invariant differential operators of degree $\le j$, then
$N(\ph)=\sum_{j=1}^n N_{D_j}(\ph)$ is a norm that makes $\CC^{j,\mu,-}(A)$ into a Banach space.
A different choice of basis will give an equivalent norm.

\begin{theorem}\label{lefschetz} (Lefschetz Formula)\\
There exists $j\in\N$ and $\mu\in\a^{*}$ such that
for any $\ph\in \CC^{j,\mu,-}(A)$ we have
$$
\sum_{\pi\in\hat G} N_\Ga(\pi)\sum_{\la\in\a^*} m_\la^\sigma(\pi) \int_{A^-}\ph(a)
a^{\la+\rho} da\=\sum_{[\ga]\in\CE(\Ga)}\ind(\ga)\,
\tr\sigma(m_\ga)\,\ph(a_\ga),
$$
where all sums and integrals converge absolutely. The inner sum on the left is always finite,
more precisely it has length $\le |W(\h,\g)|$. The left hand side is called the
\emph{global side} and the other the \emph{local side} of the Lefschetz Formula. 
Both sides of the formula give a continuous linear functional on the Banach space $\CC^{j,\mu,-}(A)$.
\end{theorem}

\prf
The Selberg trace
formula \cite{wall} says that for a compactly supported function $f$ on $G$ that is $(\dim
G+1)$-times continuously differentiable one has
$$
\sum_{\pi\in\hat G} N_\Ga(\pi)
\tr\pi(f)\=\sum_{[\ga]}\vol(\Ga_\ga\bs G_\ga)\, \CO_\ga(f),
$$
where the sum on the right hand side runs over the
conjugacy classes of $\Ga$, and $\CO_\ga(f)$ is the orbital integral
$\CO_\ga(f)\=\int_{G/G_\ga}f(x\ga x^{-1}) dx$. We need to extend the trace formula beyond
compactly supported functions.

For $d\in\N$ let $\CC^{2 d}(G)$ be the space of all functions $f$ on $G$ which are $2 d$-times
continuously differentiable and satisfy
$D f\in L^1(G)$ for every $D\in U(\g)$ with $\deg D\le 2 d$. 

\begin{lemma}
Let $d>\frac{\dim G}{2}$ and let $f\in \CC^{2 d}(G)$.
Then the trace formula is valid for $f$.
\end{lemma}

\prf In \cite{geom} Lemma 2.6 the lemma is proven in the following two cases:
\begin{enumerate}
  \item $f\ge 0$, or
  \item the geometric side of the trace formula converges with $f$ replaced by $|f|$.
\end{enumerate}
Now let $f$ be arbitrary. 
First there is $\tilde f \in \CC^{2 d}(G)$ such that $\tilde f\ge |f|$. To construct such a
function one proceeds as follows. First choose a smooth function $b:\R\ra [0,\infty)$ with
$b(x)=|x|$ for $|x|\ge 1$ and $b(x)\ge|x|$ for every
$x$. Next let $b(\eps,x)=\eps b(x/\eps)$ for $\eps>0$. Choose a function $h\in \CC^{2 d}(G)$
with $0<h(x)\le 1$ and let $\tilde f(x)=b(h(x),f(x))$. Then $\tilde f\in\CC^{2 d}(G)$ and
$\tilde f\ge |f|$. By (a) the trace formula is valid for $\tilde f$, hence the geometric
side converges for $\tilde f$, so it converges for $|f|$, so, by (b) the trace formula is
valid for $f$.
\qed

We now take $j\in \N$ and $\mu\in\a^{*}$ and keep these at our disposal. We let $\ph\in
\CC^{j,\mu,-}(A)$ and we construct a test function $f$ which has the following orbital
integrals for semisimple $y\in G$. If $y$ is not conjugate to an element of $A^- M$, then
$\CO_y(f)=0$. For
$am\in A^-M$ one has
$$
\CO_{am}(f)\= \frac{{\tr\sigma(m)}}{\det(1-am\mid \n)}\ph(a).
$$
The construction of such a function $f$ is rather straightforward.
One chooses a function $\eta\in C_c^\infty(N)$ with
$\eta\ge 0$ and $\int_N\eta(n)d n=1$. Then one sets
$$
f(k n\, am\, (k n)^{-1})\= \eta(n) {\tr\sigma(m)}
\frac{\ph(a)}{\det(1-am\mid\n)},
$$
for $am\in A^-M$, $k\in K$ and $n\in N$; further, one sets
$f(x)=0$ if $f$ is not conjugate to an element of $A^-M$.
We have to show that $f$ is well defined. This follows
from the next lemma.

\begin{lemma}
Let $am,a' m'\in A^-M$ and $x\in G$  suppose that 
$a' m'=x a m x^{-1}$. Then $a=a'$ and $x\in AM$.
\end{lemma}

\prf Every eigenvalue of $\Ad(am)$ on $\n$ is less than
$1$ in absolute value. The space $\n$ is the maximal subspace of $\g$ with this property. 
The same holds
for $am$ replaced by $a' m'$. Therefore $\Ad(x)$ must preserve $\n$ and so
$x$ lies in the normalizer of $\n$ which is $P=NMA$.
Suppose $x=n m_1 a_1$ and write $\hat m=m_1 m m_1^{-1}$. Then
$
x a m x^{-1}\= n a\hat m n^{-1}\= a\hat m ((a\hat m)^{-1} n
a\hat m) n^{-1}.
$
Since this lies in $A^-M$ and $AM\cap N=\{ 1\}$ we infer
that $(a\hat m)^{-1} n a\hat m=n$, which implies that
$n=1$. The lemma follows.
\qed

Let $d\in\N$. We will show that for $\Re(\mu)$ and $j$ sufficiently large  the function $f$
lies in
$\CC^{2 d}(G)$.

Note that the factor $1/{\det(1-am \mid \n)}$ has a pole at $am=1$ of order equal to the dimension of $\n$. 
To make this more precise let $r=\dim A$ and let $\beta_1,\dots,\beta_r$ be the simple roots of $(A,P)$.
There are $q_1,\dots,q_r\in \N$ such that $2\rho=q_1\beta_1+\cdots +q_r\beta_r$.
Define coordinates $a_1,\dots,a_r$ on $A$ by $a_j=-\beta_j(\log a)$.
Then $a$ lies in $A^-$ if an only if the coordinates $a_j$ are all $>0$.
The function on $A$ given by $a\mapsto \frac{a_1^{q_1}\cdots a_r^{q_r}}{\det(1-am |\n)}$ extends continously to the boundary of $A^-$, and the $q_j$ are minimal with this property.
So, if
we assume $j>\dim\n$, the function $f$ will on $AM$ be  $(j-1-\dim\n)$-times continuously
differentiable.
To investigate the differentiability on $G$ we need to look at the conjugation map.

Consider the map
$$
\begin{array}{c c c c}
F\ :& K\times N\times M\times A &\ra& G\\
{}& (k,n,m,a)&\mapsto&k n a m(k n)^{-1}.
\end{array}
$$
Then $f$ is a $j-1-\dim(\n)$-times continuously differentiable
function on $K\times N\times M\times A$ which factors over $F$.
To compute the order of differentiability as a function on $G$ we
have to take into count the zeroes of the differential of $F$. So
we compute the differential $F_*$ of $F$ which we view as a map on tangent spaces. Let at first $X\in \k$, then
$$
F_*(X) f(k n a m(k n)^{-1}) = \frac{d}{d t}|_{t=0} f(k
\exp(t X)n a m n^{-1}\exp(-t X)k^{-1}),
$$
which implies the equality
$$
F_*(X)_x = (\Ad(k)(\Ad(n(am)^{-1}n^{-1})-1)X)_x,
$$
when $x$
equals $k n a m(k n)^{-1}$. Similarly for $X\in\n$ we get that
$$
F_*(X)_x = (\Ad(k n)(\Ad((am)^{-1})-1)X)_x
$$
and for $X\in
\a\oplus \m$ we finally have $F_*(X)_x = (\Ad(k n)X)_x$.
From this it becomes clear that $F$, regular on $K\times
N\times M \times A^-$, may on the boundary have vanishing
differential of order $\dim(\n)+\dim(\k)$. Together we get
that $f$ is $(j-1-2\dim\n -\dim\k)$-times continuously differentiable on $G$. So we assume
$j\ge 2\dim(\n)+\dim(\k)+1$ from now on. In order to show
that $f$ goes into the trace formula for $\Re(\mu)$ and $j$ large
we fix $d>\frac{\dim G}{2}$. We have to show that 
$D f\in L^1(G)$ for any $D\in U(\g)$ of degree $\le 2 d$.
 For this we recall the map $F$ and our
computation of its differential. 
Let $\q\subset \k$ be a
complementary space to $\k_M$. 
On the regular set $F_*$ is
surjective and it becomes bijective on
$\q\oplus\n\oplus\a\oplus\m$. 
Fix $x=k n a m(k n)^{-1}$ in the
regular set and let $F_{*,x}^{-1}$ denote the inverse map of $F_{*,x}$
which maps to $\q\oplus\n\oplus\a\oplus\m$. Introducing norms on
the Lie algebras we get an operator norm for $F_{*,x}^{-1}$ and the
above calculations show that $\norm{F_{*,x}^{-1}}\le P(am)$, where $P$ is
a class function on $AM$, which, restricted to any Cartan $H=AB$
of $AM$ is a linear combination of quasi-characters. Supposing
$j$ and $\Re(\mu)$ large enough we get for $D\in U(\g)$ with
$\deg(D)\le 2 d$:
$$
|D f(k n a m(k n)^{-1})| \le \sum_{D_1} P_{D_1}(am) |D_1 f(k,n,a,m)|,
$$
where the sum runs over a finite set of $D_1\in
U(\k\oplus\n\oplus\a\oplus\m)$ of degree $\le 2 d$ and $P_{D_1}$
is a function of the type of $P$. On the right hand side we have
considered $f$ as a function on $K\times N\times A\times M$.
This discussion uses the facts that $K$ is compact, $N$ is unipotent, and $\det\left(\Ad(n(am)^{-1}n^{-1})-1\right)=\det\left(\Ad(am)^{-1}-1\right)$.
Finiteness of the sum in the inequality above follows from the Poincar\'e-Birkhoff-Witt Theorem.

It follows  that if $j,\Re(\mu)>>0$ then $(D_1 f) P P_{D_1}$
lies in $L^1(K\times N\times A\times M)$ for any $D_1$.  We have proven that the trace formula is valid for $f$.

We plug this test function into the trace formula and we claim that the result is precisely the Lefschetz formula.
To start with the geometric side recall that $\Ga$, being cocompact, only contains semisimple elements and that the orbital integrals over $f$ vanish except for $[\ga]\in\CE(\Ga)$, where we get
$$
\CO_\ga(f)\=\frac{\tr\sigma(m_\ga)}{\det(1-a_\ga m_\ga \mid\n)} \ph(a_\ga).
$$
For $[\ga]\in\CE(\Ga)$ we have $\la_\ga=\vol(G_\ga/\Ga_\ga)$.
So we see that
$$
\vol(\Ga_\ga\bs G_\ga) \CO_\ga(f)\= \ind(\ga) \tr\sigma(m_\ga)\ph(a_\ga),
$$
which means that the geometric side of the trace formula is the local side of the local side of the Lefschetz formula.

To compute the spectral side of the trace formula let $\pi\in\hat G$. 
Harish-Chandra showed that there is a locally integrable function $\Theta_\pi^G$ on $G$,
called the global character of $\pi$, such that
$
\tr\pi(h)\= \int_G h(x)\Theta_\pi^G(x)dx
$
for every $h\in C_c^\infty$. It follows that $\Theta_\pi^G$ is invariant under conjugation.
Hecht and Schmid have shown in \cite{HeSch} that for $am\in A^-M$, 
$$
\Theta_\pi^G\= \frac{\sum_{q=0}^{\dim\n} (-1)^q
\Theta_{H_q(\n,\pi_K)}^{AM}(am)}{\det(1-am\mid \n)},
$$
where  $\Theta^{AM}$ is the corresponding global character on the group $AM$.

Let $g\in L^1(G)$ be supported in the set $\{ x a m x^{-1} | am\in A^-M, x\in G\}$. Then, as a consequence of the Weyl integration formula or by direct proof one gets that
$\int_G g(x) dx$ equals
$$
\int_{KN}\int_{A^-M} g(k n am (k n)^{-1})\, |\det(1-am | \n\oplus\bar\n)|\, da d m d k d n.
$$
We apply this to $g(x)=\Theta_\pi^G (x)f(x)$ to get
\begin{eqnarray*}
\tr \pi(f)& = & \int_B\Theta_\pi^G(x) f(x) dx \\
 & = & \int_{A^-M}\Theta_\pi^G(am)\tr\sigma(m)\ph(a) |\det(1- am |\bar\n)| da d m.
\end{eqnarray*}
Using the result of Hecht and Schmid we see that this equals
$$
\int_{A^-M}\tr\sigma(m)\ph(a)\sum_{p=0}^{\dim \n} (-1)^p \Theta_{H_p(\n.\pi_K)}^{AM} (am)\frac{|\det(1-am |\bar\n)|}{ \det(1-am |\n)} da d m.
$$
For $am\in A^-M$ we have
\begin{eqnarray*}
|\det(1-am|\bar\n)|& = & (-1)^{\dim \n}\det(1-am|\bar\n) \\
& = &   (-1)^{\dim \n}a^{-2\rho}det(a^{-1}-m|\bar\n)\\
& = &   (-1)^{\dim \n}a^{-2\rho}det((am)^{-1}-1|\bar\n)\\
& = &   a^{-2\rho}det(1-(am)^{-1}|\bar\n)\\
& = &   a^{-2\rho}det(1-am|\n)
\end{eqnarray*}
So that
$$
\tr\pi(f)\=\int_{A^-M}\tr\sigma(m_\ga)\sum_{p=0}^{\dim \n} (-1)^p \Theta_{H_p(\n,\pi_K)}^{AM}
(am)a^{-2\rho} \ph(a) da d m.
$$
Lemma \ref{1.1} implies
$$
\sum_{p=0}^{\dim \n} (-1)^p \Theta_{H_p(\n,\pi_K)}^{AM}(am)a^{-2\rho}\= 
\sum_{p=0}^{\dim \n}(-1)^{p+\dim \n} \Theta_{H^p(\n,\pi_K)}^{AM}(am).
$$
And so
\eqn
\tr\pi(f) &=& \int_{A^-M}\sum_{p=0}^{\dim \n}(-1)^{p+\dim \n} \Theta_{H^p(\n,\pi_K)}^{AM}(am)\ph(a){\tr\sigma(m)} da d m\\
&=& \sum_{\la\in\a^*}\sum_{p=0}^{\dim \n} (-1)^{p+\dim \n}
\dim(H^p(\n,\pi_K)_\la\otimes\breve\sigma)^M \int_{A^-}a^{\la+\rho} \ph(a) da.
\endeqn
the convergence of the trace formula implies that for given $\la\in\a^*$ the number
$$
N_\Ga(\pi)\sum_{p=0}^{\dim \n} (-1)^{p+\dim \n}\dim(H^p(\n,\pi_K)_\la\otimes\breve\sigma)^M
$$
is nonzero only for finitely many $\pi\in\hat G$. Thus the spectral side of the trace formula gives the global side of the Lefschetz formula which therefore is proven. The continuity also follows from the proof.
\qed

\section{The Dirichlet series}\label{sec2}
Let $r=\dim A$ and for $k=1,\dots ,r$ let $\al_k$ be a positive real multiple of a
simple root of $(A,P)$ such that the
modular shift $\rho$ satisfies
$$
2\rho\= \al_1+\dots+\al_r.
$$
This defines $\al_1,\dots,\al_r$ uniquely up to order.
We fix a
Haar measure (i.e., a form
$B$) such that the subset of A,
$$
\{ a\in A \mid 0\le \al_k(\log a)\le 1,\ k=1,\dots,r\}
$$
has volume $1$.

For $a\in A$ and $k=1,\dots r$ let $l_k(a)=|\al_k(\log a)|$ and $l(a)=l_1(a)\cdots l_r(a)$.
For $s=(s_1,\dots,s_r)\in\C^r$ and $j\in\N$ define
$$
L^j(s)\=\sum_{[\ga]\in\CE(\Ga)} \ind(\ga)\, 
l(a_\ga)^{j+1}\,a_\ga^{s\cdot
\al},
$$
where $s\cdot \al=s_1\al_1+\dots +s_r\al_r$. We will show that this series converges if
$\Re(s_k)>1$ for $k=1,\dots,r$. Let
$D$ denote the differential operator
$$
D\= (-1)^r\left(\frac \partial{\partial s_1}\dots \frac \partial{\partial s_r}\right).
$$

Let $\hat G(\Ga)$ denote the set of all $\pi\in\hat G$, $\pi\ne triv$ with $N_\Ga(\pi)\ne
0$. For given $\pi\in\hat G$ let $\Lambda(\pi)$ denote the set of all $\la\in\a^*$ with
$m_{\la -\rho}(\pi)\ne 0$. Then $\Lambda(\pi)$ has at most $|W(\h,\g)|$ elements.

Let
$\la\in\a^*$. Since
$\al_1,\dots,\al_r$ is a basis of
$\a^*$ we can write
$\la=\la_1\al_1+\dots+\la_r\al_r$ for uniquely determined $\la_k\in\C$. 

Let $R_k(s)$, $k\in\N$ be a sequence of rational functions on $\C^r$. For an open set
$U\subset\C^r$ let $\N(U)$ be the set of natural numbers $k$ such that the pole-divisor of
$R_k$ does not intersect $U$.
We say that the series
$$
\sum_kR_k(s)
$$
\emph{converges weakly locally uniformly on} $\C^r$ if for every  open $U\subset\C^r$
the series 
$$
\sum_{k\in\N(U)} R_k(s)
$$
converges locally uniformly on $U$.

\begin{theorem}\label{2.4}
For $j\in\N$ large enough the series $L^j(s)$ converges locally uniformly in the set
$$
\{ s\in\C : \Re(s_k)>1,\ k=1,\dots, r\}.
$$
The function $L^j(s)$ can be written as Mittag-Leffler series,
\begin{eqnarray*}
L^j(s)&=& D^{j+1} \frac 1{(s_1-1)\cdots (s_r-1)}\\
&&+\sum_{\pi\in\hat G(\Ga)} N_\Ga(\pi)
\sum_{\la\in\Lambda(\pi)} m_{\la-\rho}(\pi) D^{j+1}\frac 1{(s_1+\la_1)\cdots (s_r+\la_r)}.
\end{eqnarray*}
The double series converges weakly locally uniformly on $\C^r$. For $\pi\ne triv$ and
$\la\in\Lambda(\pi)$ we have $\Re(\la_k)\ge -1$ for $k=1,\dots ,r$ and there is
$k\in\{1,\dots,r\}$ with $\Re(\la_k)>-1$. So in particular, the double series converges
locally uniformly on $\{\Re(s_k)>1\}$.
\end{theorem}

The proof will occupy the rest of this section.
We will show that the series $L^j(s)$ converges if the real parts $\Re(s_k)$ are
sufficiently large for $k=1,\dots,r$. Since $L^j(s)$ is a Dirichlet series with positive
coefficients, the convergence in the set $\{\Re(s_k)>1\}$ will follow, once we have
established holomorphy there. This holomorphy will in turn follow from the convergence of
the Mittag-Leffler series.

Let
$$
\a_\R^{*,+}\=\{ \la_1\al_1+\dots+\la_r\al_r\mid \la_1,\dots,\la_r >0\}
$$
be the dual positive cone.
Let $\overline{\a_\R^{*,+}}$ be the closure of $\a_\R^{*,+}$ in $\a_\R^*$.

\begin{proposition}\label{2.2}
 Let $\pi\in\hat G$,
$\la\in\a^*$ with
$m_{\la}(\pi)\ne 0$. Then $\Re(\la)$ lies in the set
$$
C\=  -3\rho+\overline{\a_\R^{*,+}}.
$$
For $\pi\in\hat G$ and $\Re(\la)=-3\rho$ we have $m_{\la}(\pi)=0$ unless $\pi$ is the
trivial representation and $\la=-3\rho$ in which case $m_{\la}(\pi)=1$. 
\end{proposition}

\prf
We introduce a partial order on $\a^*$ by

\begin{tabular}{ccl}
$\mu > \nu$ & $\Leftrightarrow$ & $\mu-\nu$ is a linear combination,\\
&& with positive
integral coefficients, of roots in $\Phi^+$.
\end{tabular}

\begin{lemma}\label{vanishing}
Let $p\in\N$, let $\pi\in\hat G$ and $\mu\in\a^*$ such that $H_p(\n,\pi_K)_\mu\ne 0$. Then
there exists $\nu\in\a^*$ with $\nu<\mu$ and $H_0(\n,\pi_K)_\nu\ne 0$.

Equivalently, if $0\le p<d=\dim(\n)$ and $H^p(\n,\pi_K)_\mu\ne 0$, then there exists
$\eta\in\a^*$ with 
$\eta <\mu$ and
$H^d(\n,\pi_K)_\eta\ne 0$.
\end{lemma}
\prf
The first assertion is a weak version of Proposition 2.32 in \cite{HeSch} and the second
follows from the first and Lemma \ref{1.1}.
\qed

Theorem 4.25 of \cite{HeSch} states that the set of $\nu\in\a^*$ with $H_0(\n,\pi_K)_\nu\ne
0$ which are minimal with respect to $<$ are precisely the leading exponents of the
asymptotic of matrix coefficients of $\pi$. Since $\pi$ is unitary, its matrix coefficients
are bounded and thus its leading coefficients in the normalization of Theorem
4.16 of \cite{HeSch} all lie in the set $-\rho +\overline{\a_\R^{*,+}}$. Now Lemma \ref{vanishing} and Lemma \ref{1.1} imply the first statement of Proposition \ref{2.2}.

It remains to consider the case $\Re(\la)=-3\rho$. So let $\la=-3\rho+i\mu$ for
$\mu\in\a_\R^*$. Using the definition of Lie algebra homology it is easy to show that
$m_{-3\rho}(triv)=1$.  Let
$\pi\in\hat G$ and assume that
$H^\bullet(\n,\pi_K)^M_{-3\rho+i\mu}\ne 0$. The claim will follow, if we show that this
implies 
$\pi= triv$.

Since $H^\bullet(\n,\pi_K)\cong H_\bullet(\n,\pi_K)\otimes \bigwedge^{top}\n$ we find that our
condition is equivalent to 
$H_\bullet(\n,\pi_K)_{-\rho+i\mu}\ne 0$.

Let $\xi$ be an irreducible representation of $M$ and let $\nu\in\a^*$. let $\pi_{\xi,\nu}$ be the
induced principal series representation.  Recall that if $\Re(\nu)\in\a_\R^{*,+}$, then
$\pi_{\xi,\nu}$ has a unique irreducible quotient, its Langlands quotient. The representation
dual to $\pi_{\xi,\nu}$ is $\pi_{\breve\xi,-\nu}$, where $\breve\xi$ is the dual to $\xi$. Therefore
$\pi_{\breve\xi,-\nu}$ has a unique irreducible subrepresentation. Replacing $\xi$ by
$\breve\xi$ and $\nu$ by $-\nu$ it follows that for $\Re(\nu)\in-\a_\R^{*,+}$ the principal
series representation $\pi_{\xi,\nu}$ has a unique irreducible subrepresentation, namely the
dual of the Langlands quotient of $\pi_{\breve\xi,-\nu}$.

 Frobenius reciprocity
says that
$$
\Hom_G(\pi,\pi_{\xi,\nu})\ \cong\ \Hom_{AM}(H_0(\n,\pi_K),\xi\otimes(\nu+\rho)).
$$
So in particular,
$$
\Hom_G(\pi,\pi_{1,\nu})\ \cong\ \Hom_{A}(H_0(\n,\pi_K)^M,\nu+\rho).
$$

If $H_p(\n,\pi_K)^M_{-\rho+i\mu}\ne 0$ for some $p>0$ then by Schmid's vanishing result it
follows that $\pi$ must have a leading coefficient with real part $<-3\rho$ which we already
have ruled out. So it follows that $H_0(\n,\pi_K)^M_{-\rho+i\mu}\ne 0$ and hence
$\Hom_G(\pi,\pi_{1,-\rho+i\mu})\ne 0$. None of the representations $\pi_{1,-\rho+i\mu}$,
however, has a unitary subrepresentation, except for $\mu=0$.
 
The representation
$\pi_{1,-\rho}$ has a unique irreducible subrepresentation which is the trivial representation. 
This implies $\pi=triv$.
\qed 

For $a\in A$ set
$$
\ph(a)\= l(a)^{j+1}\ a^{s\cdot \al}.
$$
For $\Re(s_k)>>0$, $k=1,\dots,r$ the Lefschetz formula is valid for this test function. The local
side of the Lefschetz formula for $\sigma= triv$ equals
$$
\sum_{[\ga]\in\CE(\Ga)} \ind(\ga)\, l(a_\ga)^{j+1}\ a_\ga^{s\cdot
\al}\=  L^j(s).
$$
The convergence assertion in the Lefschetz formula implies that the series  converges
absolutely  if
$\Re(s_k)$ is sufficiently large for every
$k=1,\dots r$.  Since $L^j(s)$ is a Dirichlet series with positive coefficients it will
converge locally uniformly for $s$ in some open set. We will show that it extends to a
holomorphic function in  the set
$\Re(s_k)>1$,
$k=1,\dots ,r$. Again, since it is a Dirichlet series with positive coefficients it must
therefore converge in that region.

 With our given test
function and the Haar measure chosen we compute
\eqn
\int_{A^-}\ph(a) a^\la da &=&(-1)^{r(j+1)} \int_{A^-} (\al_1(\log a)\cdots\al_r(\log
a))^{j+1} a^{s\cdot
\al+\la} da\\
\endeqn
\vspace{-30pt}
\eqn
&=& (-1)^{r(j+1)} \int_0^\infty\dots\int_0^\infty \left( {t_1}\cdots t_r\right)^{j+1}
e^{-((s_1+\la_1)t_1+\dots+(s_r+\la_r)t_r)} d t_1\dots d t_r\\
&=& D^{j+1} 
\int_0^\infty\dots\int_0^\infty  e^{-((s_1+\la_1)t_1+\dots+(s_r+\la_r)t_r)} d t_1\dots d t_r\\
&=& D^{j+1}  \frac
1{(s_1+\la_1)\dots(s_r+\la_r)}.
\endeqn
Writing $m_\la=m_\la^{triv}$ and performing a $\rho$-shift we see that the Lefschetz formula
gives
\begin{eqnarray*}
L^j(s)&=&\sum_{\pi\in\hat G}N_\Ga(\pi)\sum_{\la\in\a^*} m_{\la-\rho}(\pi)\, D^{j+1} \frac
1{(s_1+\la_1)\cdots(s_r+\la_r)}\\
&=& \sum_{\pi\in\hat G}N_\Ga(\pi)\sum_{\la\in\a^*} m_{\la-\rho}(\pi)\,  \frac
{((j+1)!)^r}{(s_1+\la_1)^{j+2}\cdots(s_r+\la_r)^{j+2}}
\end{eqnarray*}
for $\Re(s_k)>>0$.
For every $\pi\in\hat G$ we fix a representative $\La_\pi\in(\a +\t)^*$ of the
infinitesimal character of $\pi$. According to Lemma \ref{1.1}, if
$m_{\la-\rho}(\pi)\ne 0$, then $\la-w\La_\pi|_\a-\rho$ for some $w\in W(\h,\g)$. By abuse of
notation we will write $w\La_\pi$ instead of $w\La_\pi|_\a$. Hence we get
$$
L^j(s)\=\sum_{\pi\in\hat G}N_\Ga(\pi)\sum_{w\in W(\h,\g)} m_{w\La_\pi-2\rho}(\pi)\, D^{j+1}
\frac 1{(s_1+\la_1)\cdots(s_r+\la_r)}.
$$ 

For $\la\in\a^*$ let $\norm\la$ be the norm given by the form $B$ as explained in the
beginning of section \ref{sec1}.

\begin{proposition}\label{2.5}
There are $m\in\N$, $C>0$ such that for every $\pi\in\hat G$ and every $\la\in\a^*$ one has
$$
|m_{\la-\rho}(\pi)|\ \le\ C(1+\norm\la)^m.
$$
\end{proposition}

\prf
Recall from the first section that for $a\in A^-$,
\begin{eqnarray*}
\sum_{\la\in\a}m_{\la-\rho}(\pi)\, a^\la&=& \int_M\sum_{p=0}^{\dim
\n}(-1)^p\Theta_{H^p(\n,\pi_K)}^{AM}(am)\, dm\\
&=& (-1)^{\dim\n}a^{-2\rho} \int_M\sum_p(-1)^p\Theta_{H_p(\n,\pi_K)}^{AM}(am)\, dm\\
&=& (-1)^{\dim\n}a^{-2\rho} \int_M\det(1-am|\n)\,\Theta_\pi^G(am)\, dm
\end{eqnarray*}
Choose a set of positive roots $\phi^+(\t,\m)\subset\phi(\t,\m)$ compatible to the choice
of $P$. Let $\rho_M=\frac 12\sum_{\al\in \phi^+(\t,\m)}\al$. For $t\in T$ set
$$
D_T(t) \df t^{\rho_M} \prod_{\al\in \phi^+(\t,\m)} (1-t^{-\al}).
$$
This is the Weyl denominator for $T$. The Weyl integral formula gives
$$
\sum_{\la\in\a^*}m_{\la-\rho}(\pi)\, a^\la\= \int_{T^reg} (-1)^{\dim\n}a^{-2\rho}
\det(1-am|\n)|D_T(t)|^2\Theta_\pi^G(at)\, dt.
$$
The function $(-1)^{\dim\n}a^{2\rho}\det(1-am|\n)D_T(t)$ equals the Weyl denominator for
$H=AT$. By Theorems 10.35 and 10.48 of \cite{Knapp} there are constants $c_w$, $w\in
W(\h,\g)$ such that
$$
(-1)^{\dim\n}a^{-2\rho}
\det(1-am|\n)D_T(t)\Theta_\pi^G(at)\= \sum_{w\in W(\h\g)}c_w\, (at)^{w\La_\pi}.
$$
We thus have proved the following Lemma.
\begin{lemma}\label{2.6}
 For $a\in A^-$,
$$
\sum_{\la\in\a^*}m_{\la-\rho}(\pi)\, a^\la\= \sum_{w\in
W(\h,\g)}c_w\,a^{w\La_\pi}\int_{T^{reg}}
t^{w\La_\pi-\rho_M}\prod_{\al\in\phi^+(\t,\m)}(1-t^\al)\, dt. 
$$
\end{lemma}

Proposition \ref{2.5} will follow from explicit formulae for the global character
$\Theta_\pi^G$ (see below) which give bounds on the $c_w$.
Another remarkable consequence of Lemma \ref{2.6} is the fact that
there is a finite set $E\subset\t^*$ such that whenever $m_{\la-\rho}(\la)\ne 0$ for some
$\la\in\a^*$ it follows $\La_\pi|_\t\in E$. Hence Proposition \ref{2.5} will follow from the
estimate
$$
|m_{\la-\rho}(\pi)|\ \le\ C(1+\norm{\La_\pi})^m.
$$

In \cite{HC-DSII} Harish-Chandra gives an explicit formula for characters of discrete
series representations which imply the sharper estimate $|m_{\la-\rho}(\pi)|\le C$ for
the discrete series representations. From Harish-Chandra's paper a similar formula can be
deduced for limit of discrete series representations. Alternatively, one can use Zuckerman
tensoring (Prop. 10.44 of \cite{Knapp}) to deduce the estimate for limits of discrete
series representations.
Next, if $\pi=\pi_{\sigma,\nu}$ is induced from some parabolic $P_1=M_1A_1N_1$, then the
character of $\pi$ can be computed from the character of $\sigma$ and $\nu$, see formula
(10.27) in \cite{Knapp}. From this it follows that the claim holds for
standard representations, i.e.
admissible representations which are induced from discrete series or limit of discrete
series representations. 

\begin{lemma}\label{2.7}
There are natural numbers $n,m$ and a constant $d>0$ such that for every $\pi\in\hat G$
there are standard representations $\pi_1,\dots,\pi_n$ and integers $c_1,\dots,c_n$ with
$$
\Theta_\pi\= \sum_{k=1}^n c_k\,\Theta_{\pi_k}
$$
and $|c_k|\le d(1+\norm{\La_\pi}^m)$ for $j=1,\dots, n$.
\end{lemma}

\prf By the Langlands classification, $\pi$ is a quotient $p/q$ of a standard representation
$p$. Then $\pi$ and $p$ and thus every constituent of $q$ share the same infinitesimal
character $\La_\pi$. The $A$-parameter of every constituent of $q$ is smaller than the
one of $\pi$, so one can perform an induction on this $A$-parameter as suggested by
(\cite{Knapp}, Problem 10.2) to get the linear combination as above.
The length of this linear combination is globally bounded because
the combination only contains characters of representations that share a common infinitesimal
character and the number of such is globally bounded. The coefficients are bounded by the
multiplicity that a constituent $\pi$ of an induced representation
$\pi_{\sigma,\nu}$ might have, where $\sigma$ is a discrete series or limit of discrete series representation. 
Now
assume that $\pi$ contains a $K$-type $\tau$. Then the multiplicity of $\pi$ in $\pi_{\sigma,\nu}$ is bounded by the
multiplicity of $\tau$ in $\pi_{\sigma,\nu}$ which in turn is bounded by $\dim\tau$. It remains to show that every
irreducible admissible $\pi$ contains a minimal $K$-type whose dimension is bounded by a power of
$\norm{\Lambda_\pi}$. By Weyl's character formula this is equivalent to say that the norm of the infinitesimal character
of $\tau$ is bounded by a power of $\norm{\Lambda_\pi}$. This indeed follows from the minimal $K$-type formula
(\cite{Knapp}, Theorem 15.1) together with (\cite{Knapp}, Theorem 15.10). So we have shown that the 
$m_{\la-\rho}(\pi)$ grow
at most like a power of $\norm{\La_\pi}$ and Lemma \ref{2.7} follows. By the above, this
also implies Proposition \ref{2.5}.
\qed

It remains to deduce Theorem \ref{2.4}. Since the coefficients $m_{\la-\rho}(\pi)$ grow at
most like a power of $\norm{\La_\pi}$, the convergence assertion in Theorem \ref{2.4} will
be implied by the following lemma.

\begin{lemma}\label{2.8}
Let $S$ denote the set of all pairs $(\pi,\la)\in\hat G\times\a^*$ such that
$m_{\la-\rho}(\pi)\ne 0$. There is $m_1\in\N$ such that
$$
\sum_{(\pi,\la)\in S} \frac{N_\Ga(\pi)}{(1+\norm{\la})^{m_1}}\ <\ \infty.
$$
\end{lemma}

\prf
By the remark following Lemma \ref{2.6} it suffices to show that 
there is $m\in\N$ such that
$$
\sum_{(\pi,\la)\in S} \frac{N_\Ga(\pi)}{(1+\norm{\La_\pi})^{m_1}}\ <\ \infty.
$$
Let $\pi\in\hat G$. The restriction of $\pi$ to the maximal compact subgroup $K$ decomposes
into finite dimensional isotypes
$$
\pi|_K\=\bigoplus_{\tau\in\hat K} \pi(\tau).
$$
Let $C_K$ be the Casimir operator of $K$ and let
$$
\Delta_G\df -C+2C_K.
$$
Then $\Delta_G$ is the Laplacian on $G$ given by the left invariant metric which at the
point $e\in G$ is given by
$\sp{.}{.}=-B(.,\theta(.))$. Since $\Delta_G$ is left invariant it induces an operator on
$\Ga\bs G$ denoted by the same letter. This operator is
$\ge 0$ and elliptic, so there is a natural number $k$ such that $(1+\Delta_G)^{-k}$ is of
trace class on $L^2(\Ga\bs G)$.
Hence
\begin{eqnarray*}
\infty &>& \tr(1+\Delta_G)^{-k}\\
&=& \sum_{\pi\in\hat G} N_\Ga(\pi) \sum_{\tau\in\hat K}(1-\pi(C)+2\tau(C_K))^{-k}\,
\dim\pi(\tau)\\
&\ge& \sum_{\pi\in\hat G} \frac{N_\Ga(\pi)}{(1- \pi(C)+2\tau_\pi(C_K))^k},
\end{eqnarray*}
where  for each $\pi\in\hat G$ we fix a minimal $K$-type $\tau_\pi$.
Since the infinitesimal character of the minimal $K$-type grows like the infinitesimal
character of $\pi$ the Lemma follows.
\qed

Finally, to prove Theorem \ref{2.4}, let $U\in\C^r$ be open. Let $S(U)$ be the set of all
pairs $(\pi,\la)\in\hat G\times\a^*$ such that $m_{\la-\rho}(\pi)\ne 0$ and the pole divisor
of
$$
\frac 1{(s_1+\la_1)\cdots (s_r+\la_r)}
$$
does not intersect $U$.
Let $V\subset U$ be a compact subset. We have to show that for some $j\in\N$ which does not
depend on $U$ or $V$,
$$
\sup_{s\in V}\sum_{(\pi,\la)\in S(U)} \left| \frac{N_\Ga(\pi) \,
m_{\la-\rho}(\pi)}{(s_1+\la_1)^{j+2}\cdots (s_r+\la_r)^{j+2}}\right|\ <\ \infty.
$$
Let $m$ be as in Lemma \ref{2.5} and $m_1$ as in Lemma \ref{2.8}. Then let $j\ge m+m_1-2$.
Since $V\subset U$ and $V$ is compact there is $\eps>0$ such that $s\in V$ and
$(\pi,\la)\in S(U)$ implies $|s_k+\la_k|\ge\eps$ for every $k=1,\dots,r$. Hence there is
$c>0$ such that for every $s\in V$ and every $(\pi,\la)\in S(U)$,
$$
|(s_1+\la_1)\cdots (s_r+\la_r)|\ \ge\ c(1+\norm\la).
$$
This implies,
\begin{eqnarray*}
\left|\frac{m_{\la-\rho}(\pi)}{(s_1+\la_1)^{j+2}\cdots (s_r+\la_r)^{j+2}}\right|&\le& \frac
1{c^{j+2}}\ \frac{|m_{\la-\rho}(\pi)|}{(1+\norm\la)^{j+2}}\\
&\le& \frac C{c^{j+2}}\ \frac{1}{(1+\norm\la)^{j+2-m}}\\
&\le& \frac C{c^{j+2}}\ \frac{1}{(1+\norm\la)^{m_1}}.
\end{eqnarray*}
The claim now follows from Lemma \ref{2.8}. The proof of Theorem \ref{2.4} is finished.
\qed

\section{The prime geodesic theorem}\label{sec3}

We now give the main result of the paper.

\begin{theorem} (Prime Geodesic Theorem)\\
For $T_1,\dots,T_r >0$ let
$$
\Psi(T_1,\dots,T_r)\=\sum_{\stackrel{\stackrel{[\ga]\in\CE(\Ga)}{}}{\stackrel{}{a_\ga^{-\al_k}\le
T_k,\ k=1,\dots,r}}}
\la_\ga.
$$
Then, as $T_k\to\infty$ for $k=1,\dots,r$ we have
$$
\Psi(T_1,\dots,T_r)\ \sim\ T_1\cdots T_r.
$$
\end{theorem}

The proof of the theorem will occupy the rest of the section.

For $x_1,\dots,x_r>0$ let
$$
A(x)\= A(x_1,\dots,x_r)\= \sum_{\stackrel{[\ga]\in\CE(\Ga)}{l_k(a_\ga)\le x_k}}
l(a_\ga)^{j+1}\,\ind(\ga).
$$
Let $\R_+$ be the set of positive real numbers. Then
$
L^j(s)\= \int_{\R_+^r} A(x) e^{-s\cdot x} dx.
$

We need a higher dimensional analogue of the Wiener-Ikehara Theorem. For this we introduce a
partial order on $\R^r$. We define
$$
(x_1,\dots,x_r)\ge (y_1,\dots,y_r)\ \Leftrightarrow\ x_1\ge y_1\ {\rm and}\dots {\rm and}\ x_r\ge y_r.
$$
A real valued function $F$ on a subset of $\R^r$ is called \emph{monotonic} if $x\ge y$ implies
$F(x)\ge F(y)$. The partial order also gives sense to the assertion that a sequence $x^n$ tends to
$+\infty$. It does so if all its components $x_k^n$ tend to $+\infty$ separately.

\begin{theorem}
(Higher dimensional version of the Wiener-Ikehara Theorem)\\
Let $A\ge 0$ be a monotonic measurable function on $\R_+^r$. Suppose that the integral
$
L(s)\= \int_{\R_+^r} A(x) e^{-s\cdot x} dx
$
converges for $\Re(s_k)>1$, $k=1,\dots,r$.
Suppose further that there is a countable set $I$ and for each $i\in I$, there is $\theta_{i}\in\C^r$ such that
$\Re(\theta_{i,k})\le 1$ for $k=1,\dots r$. Assume that for each $i\in I$ there is $k\in\{ 1,\dots,r\}$ such that
$\Re(\theta_{i,k})< 1$ and that there are integers $c_{i}$ such that the function
$$
L(s)-D^{j+1}\frac 1{(s_1-1)\cdots (s_r-1)}
-\sum_{i\in I} c_{i}\, D^{j+1}\frac
1{(s_1-\theta_{i,1})\cdots (s_r-\theta_{i,r})}
$$
extends to a  holomorphic function on the set of all $s\in\C^r$ with 
$ \Re(s_k)\ge 1$ for $k=1,\dots,r$. Here we assume that the sum converges weakly locally
uniformly absolutely on $\C^r$. Under these circumstances we can
conclude 
$$
\lim_{x\to +\infty} A(x)\, (x_1\cdots x_r)^{-(j+1)}e^{-(x_1+\dots+x_r)}\= 1.
$$
\end{theorem}

The proof of this theorem is a fairly straightforward application of methods from analytic
number theory. We include it for the convenience of the reader.

\prf
Let $A$ as in the theorem and let 
$$
B(x)=A(x)\, (x_1\cdots x_r)^{-(j+1)}e^{-(x_1+\dots+x_r)}.
$$
To simplify notations we will frequently identify a complex number $z$ with the vector
$(z,\dots,z)\in\C^r$. Since
$$
\frac 1{(s_1-1)\cdots(s_r-1)}\=\int_{\R_+^r}e^{-(s-1)\cdot x}\, dx
$$
we get
$$
D^{j+1}\frac 1{(s_1-1)\cdots(s_r-1)}\=\int_{\R_+^r}(x_1\cdots x_r)^{j+1}\, e^{-(s-1)\cdot
x}\, dx.
$$
Likewise,
$$
D^{j+1}\frac 1{(s_1-\theta_{i,1})\cdots(s_r-\theta_{i,r})}\=\int_{\R_+^r}(x_1\cdots
x_r)^{j+1}\, e^{-(s-\theta_{i})\cdot x}\, dx.
$$
Let $f$ be a smooth function of compact support
of $\R$ which is real valued and even. Then its Fourier transform $\hat f$ will also be real
valued. We further assume $f$ to be of the form $f=f_1*f_1$ for some $f_1$. Then $\hat
f=(\hat f_1)^2$ is positive on the reals. Let $I(f)$ be the set of $i\in I$ such that 
$\Im(\theta_{i})\in (\supp f)^r$.  Then the function
\begin{eqnarray*}
g(s)&=&L(s)-D^{j+1}\frac 1{(s_1-1)\cdots (s_r-1)}\\
&& -\sum_{i\in I(f)} c_{i}\, D^{j+1}\frac
1{(s_1-\theta_{i,1})\cdots (s_r-\theta_{i,r})}
\end{eqnarray*}
extends to an analytic function on $\{ \Re(s_k)\ge 1,\ \Im(s_k)\in\supp f\}$. It follows
\begin{eqnarray*}
g(s)&=&\int_{R_+^r} (B(x)-1)\, (x_1\cdots x_r)^{j+1}\, e^{-(s-1)\cdot x}\, dx\\
&& -\sum_{i\in I(f)} c_{i}\int_{\R_+^r}(x_1\cdots x_r)^{j+1}\,
e^{-(s-\theta_{i})\cdot x}\, dx
\end{eqnarray*}

Let $\eps>0$. For
$y\in\R$ the integral
$$
\int_{\R^r} g(1+\eps+it) f(t_1)\cdots f(t_r)\, e^{iy\cdot t}\, dt_1\dots dt_r
$$
equals
$$
 \int_{\R^r}  f(t_1)\cdots f(t_r)\, e^{iy\cdot t}\, \int_{\R_+^r}(B(x)-1)\, (x_1\cdots
x_r)^{j+1}\, e^{-(\eps+it)\cdot x}\, dx\, dt
$$ $$
-\sum_{i\in I(f)} c_{i} \int_{\R^r}  f(t_1)\cdots f(t_r)\, e^{iy\cdot t}\, 
\int_{\R_+^r}(x_1\cdots x_r)^{j+1}\,
e^{-(1-\theta_{i})\cdot x}\, e^{-\eps\cdot x}dx\, dt.
$$
We want to interchange the order of integration. This only causes a problem for the summand involving $B(x)$. To
justify the interchange, note that by the monotonicity of
$A$ we have for real $s$, and $x\in\R_+^r$,
$$
L(s)\=\int_{\R_+^r} A(u)\, e^{-s\cdot u}\, du\ \ge\ A(x)\int_{x+\R_+^r}\int_{x+\R_+^r} e^{-s\cdot
u}\, du\=\frac{A(x)\, e^{-x\cdot s}}{s_1\cdots s_r}.
$$
In other words, $A(x)\le s_1\cdots s_r L(s) e^{x\cdot s}$. Therefore $A(x)=O(e^{-s\cdot x})$ for
every $s_1,\dots,s_r>1$ which implies $A(x)=o(e^{-s\cdot x})$ for
every $s_1,\dots,s_r>1$. So for $\delta>0$, $B(x)\,(x_1\cdots x_r)^{j+1}\, e^{-\delta\cdot
x}=A(x)\, e^{-(1+\delta)\cdot x}=o(1)$ for every $\delta>0$. This implies that the integral
$$
\int_{\R_+^r}(B(x)-1)\, (x_1\cdots x_r)^{j+1}\, e^{-(\eps+it)\cdot x}\, dx
$$
converges locally uniformly in $t$. So we can interchange the order of
integration to obtain that
$$
\int_{\R^r}g(1+\eps+it)\, f(t_1)\cdots f(t_r)\, e^{iy\cdot t}\, dt
$$
equals
$$
\int_{\R_+^r}(B(x)-1)\, (x_1\cdots x_r)^{j+1}\, e^{-\eps\cdot x} \hat f(y_1-x_1)\cdots\hat
f(y_r-x_r)\, dx
$$ $$
-\sum_{i\in I(f)} c_{i} \int_{\R_+^r}  (x_1\cdots x_r)^{j+1}e^{-(1-\theta_{i})\cdot
x}\hat f(y_1-x_1)\cdots
\hat f(y_r-x_r)\, \, 
e^{-\eps\cdot x}dx.
$$
Since $g(s)$ is analytic in the set 
$\Re(s_1),\dots,\Re(s_r)\ge 1$ we can let $\eps\ra 0$ to obtain that
$$
\int_{\R^r}g(t)\, f(t_1)\cdots f(t_r)\, e^{iy\cdot t}\ dt
$$
equals
$$
\lim_{\eps\ra 0}\int_{\R_+^r}(B(x)-1)\, (x_1\cdots x_r)^{j+1}\, e^{-\eps\cdot x}\,
\hat f(y_1-x_1)\cdots \hat f(y_r-x_r)\, dx.
$$
$$
-\sum_{i\in I(f)} c_{i} \int_{\R_+^r}  (x_1\cdots x_r)^{j+1}e^{-(1-\theta_{i})\cdot
x}\hat f(y_1-x_1)\cdots
\hat f(y_r-x_r)\, \, 
e^{-\eps\cdot x}dx.
$$
Since $\hat f$ is rapidly decreasing the limit for $\eps\ra 0$ of 
$$
\int_{\R_+^r}(x_1\cdots x_r)^{j+1} \,
e^{-\eps\cdot x}\, \hat f(y_1-x_1)\cdots\hat f(y_r-x_r)\, dx
$$
exists and equals
$$
\int_{\R_+^r}(x_1\cdots x_r)^{j+1} \,
 \hat f(y_1-x_1)\cdots\hat f(y_r-x_r)\, dx.
$$
For the sum over $I(f)$ recall that the imaginary parts of the $\theta$'s are in a compact set,
therefore the real parts must tend to $-\infty$ and so the convergence is uniform in $\eps$, i.e.,
the limit can be interchanged with the summation. Hence also the limit
$$
\lim_{\eps\ra 0}\int_{\R_+^r}B(x)\, (x_1\cdots x_r)^{j+1}\, e^{-\eps\cdot x}\,
\hat f(y_1-x_1)\cdots \hat f(y_r-x_r)\, dx
$$
exists. Since we assume $\hat f\ge 0$ the integrand is nonnegative and monotonically
increasing as $\eps\ra 0$. Therefore the limit may be drawn under the integral sign. We
conclude that
$$
\int_{\R^r}g(t)\, f(t_1)\cdots f(t_r)\, e^{iy\cdot t}\ dt
$$
equals
$$
\int_{\R_+^r}(B(x)-1)\, (x_1\cdots x_r)^{j+1}
\hat f(y_1-x_1)\cdots \hat f(y_r-x_r)\, dx.
$$
$$
-\sum_{i\in I(f)} c_{i} \int_{\R_+^r}  (x_1\cdots x_r)^{j+1}e^{-(1-\theta_{i})\cdot
x}\hat f(y_1-x_1)\cdots
\hat f(y_r-x_r)\, dx.
$$
By the Riemann-Lebesgue Lemma this tends to zero as $y\ra\infty$. 
For $y>>0$ we estimate
\begin{eqnarray*}
&&\int_{\R_+^r} (x_1\cdots x_r)^{j+1}\, e^{(\theta-1)\cdot x}\, \hat f(y_1-x_1)\cdots \hat
f(y_r-x_r)\, dx\\
&\le & \int_{\R^r} (x_1\cdots x_r)^{j+1}\, e^{(\theta-1)\cdot x}\, \hat f(y_1-x_1)\cdots
\hat f(y_r-x_r)\, dx\\
&=& \int_{\R^r} ((x_1+y_1)\cdots (x_r+y_r))^{j+1}\, e^{(\theta-1)\cdot (x+y)}\, \hat
f(-x_1)\cdots
\hat f(-x_r)\, dx\\
&\le& (const)\ (y_1\cdots y_r)^{j+1} \, e^{(\theta-1)\cdot y}.
\end{eqnarray*}
This implies that the sum over $I(f)$ tends to zero as $y\ra\infty$.

Therefore
$$
\lim_{y\ra\infty}\int_{\R_+^r}(B(x)-1)\, (x_1\cdots x_r)^{j+1}
\hat f(y_1-x_1)\cdots \hat f(y_r-x_r)\, dx\= 0.
$$

\begin{lemma}
For every $k=0,1,2,\dots$,
$$
\lim_{y\ra\infty}\frac 1{y^k}\int_0^\infty x^k\,\hat f(y-x)\, dx\= 2\pi f(0).
$$
\end{lemma}

\prf
Start with $k=0$. Then
\begin{eqnarray*}
\lim_{y\ra\infty}\int_0^\infty \hat f(y-x)\, dx &=& \lim_{y\ra\infty}\int_{-\infty}^\infty \hat f(y-x)\, dx\\
&=& \int_{-\infty}^\infty \hat f(-x)\, dx\= 2\pi f(0)
\end{eqnarray*}
For $k\mapsto k+1$ consider
\begin{eqnarray*}
\lim_{y\ra\infty}\frac 1{y^{k+1}} \int_0^\infty x^{k+1}\,\hat f(y-x)\, dx &=& \lim_{y\ra\infty}\frac
1{y^{k+1}} \int_{-\infty}^\infty x^{k+1}\,\hat f(y-x)\, dx\\
&=& \lim_{y\ra\infty}\frac 1{y^{k+1}} \int_{-\infty}^\infty (x+y)^{k+1}\,\hat f(-x)\, dx\\
&=& 2\pi f(0).
\end{eqnarray*}
\qed

This lemma implies that
$$
\lim_{y\ra\infty}\int_{\R_+^r}B(x)\,\left(\frac{x_1\cdots x_r}{y_1\cdots
y_r}\right)^{j+1}\, \hat f(y_1-x_1)\cdots \hat f(y_r-x_r)\, dx\= (2\pi f(0))^r.
$$
Let $S>0$. Since $A(x)$ is monotonic we have
$A(y-S)\le A(x)\le A(y+S)$ whenever $y-S\le x\le y+S$. In that range we then have
$$
B(y-S)((y_1-S)\cdots (y_r-S))^{j+1} e^{(y-S)\cdot 1}
 \le\ B(x)(x_1\cdots x_r)^{j+1}\, e^{x\cdot 1}\hspace{30pt}
$$ $$
\hspace{50pt}\le\ B(y+S)((y_1+S)\cdots (y_r+S))^{j+1} e^{(y+S)\cdot 1}.
$$
The first inequality implies
\begin{eqnarray*}
B(x)(x_1\cdots x_r)^{j+1} &\ge& B(y-S)((y_1-S)\cdots(y_r-S))^{j+1}\, e^{(y-x-S)\cdot 1}\\
&\ge& B(y-S)((y_1-S)\cdots(y_r-S))^{j+1}\, e^{-2S\cdot 1}.
\end{eqnarray*}
So for $y\ge S$,
$$
e^{-2rS}\, B({y-S})\left(\frac{(y_1-S)\cdots (y_r-S)}{y_1\cdots
y_r}\right)^{j+1}\int_{y-S}^{y+S}\hat f(y_1-x_1)\cdots \hat f(y_r-x_r)\, dx
$$
\begin{eqnarray*}
&\le& \int_{y-S}^{y+S}B(x)\left(\frac{x_1\cdots x_r}{y_1\cdots y_r}\right)^{j+1}\hat
f(y_1-x_1)\cdots \hat f(y_r-x_r)\, dx\\
&\le& \int_{\R_+^r}B(x)\left(\frac{x_1\cdots x_r}{y_1\cdots y_r}\right)^{j+1}\hat
f(y_1-x_1)\cdots \hat f(y_r-x_r)\, dx.
\end{eqnarray*}
This implies
$$
\limsup_{y\ra\infty}B(y)\ \le\ e^{2rS}\frac{(2\pi f(0))^2}{\int_{-S}^S\hat f(x_1)\cdots \hat f(x_r)\, dx}.
$$
We vary $f$ so that $\hat f$ is small outside $[-S,S]$. In this way we get
$$
\limsup_{y\ra\infty}B(y)\ \le\  e^{2rS}.
$$
Since this is true for any $S>0$ it follows
$$
\limsup_{y\ra\infty}B(y)\ \le\  1.
$$
The inequality $\liminf_{y\ra\infty}B(y)\ge 1$ is obtained in a similar fashion.
The Wiener-Ikehara Theorem is proven.
\qed

Let
$$
\phi(T)\=\phi(T_1,\dots,T_r) \df \sum_{\stackrel{[\ga]\in\CE(\Ga)}{a_\ga^{-\al_k}\le T_k}}
\frac{\la_\ga}{\det(1-a_\ga m_\ga\mid\n)}.
$$

\begin{lemma}
As $T_1,\dots T_r\ra\infty$ we have
$$
\phi(T_1,\dots T_r)\ \sim\ T_1\cdots T_r.
$$
\end{lemma}

\prf
Let  
$$
\phi_j(T)\df \sum_{\stackrel{[\ga]\in\CE(\Ga)}{a_\ga^{-\al_k}\le T_k}}
\frac{\la_\ga\left( l_1(a_\ga)\cdots l_r(a_\ga)\right)^{j+1}}{\det(1-a_\ga m_\ga\mid\n)}.
$$
By the Wiener -Ikehara theorem  we know that 
$$
\frac{\phi_j(T)}{T_1\cdots T_r (\log T_1)^{j+1}\cdots (\log T_r)^{j+1}}
$$
tends to $1$ as $T\ra\infty$. Obviously
$$
\frac{\phi_j(T)}{(\log T_1)^{j+1}\cdots (\log T_r)^{j+1}}\ \le\ \phi(T),
$$
which implies
$$
\liminf_{T\ra\infty} \frac{\phi(T)}{T_1\cdots T_r}\ \ge\ 1.
$$
Let $0<\mu<1$. Then
\begin{eqnarray*}
\phi_j(T) &\ge& \sum_{\stackrel{[\ga]\in\CE(\Ga)}{T_k^\mu < a_\ga^{-\al_k}\le T_k}}
{\ind(\ga)\left( l_1(a_\ga)\cdots l_r(a_\ga)\right)^{j+1}}\\
&\ge& \mu^{r(j+1)}\left((\log T_1)\cdots (\log
T_r)\right)^{j+1}\sum_{\stackrel{[\ga]\in\CE(\Ga)}{T_k^\mu < a_\ga^{-\al_k}\le T_k}}
{\ind(\ga)}\\ &=& \mu^{r(j+1)}\left((\log T_1)\cdots (\log T_r)\right)^{j+1}
(\phi(T)-\phi(T^\mu))
\end{eqnarray*}
From this we get
$$
\frac{\phi(T)}{T_1\cdots T_r}\ \le\ 
$$ $$
\le\ \mu^{-r(j+1)}\frac{\phi_j(T)}{T_1\cdots T_r (\log
T_1)^{j+1}\cdots (\log T_r)^{j+1}} + \frac{\phi(T^\mu)}{T_1^\mu\cdots T_r^\mu (T_1\cdots
T_r)^{1-\mu}}
$$
Assume first that $\phi(T)/T_1\cdots T_r$ tends to infinity as $T\ra\infty$. Then there is a
sequence ${_nT}$ in $\R_+^r$, tending to infinity such that $\phi({_nT})/{_nT}_1\cdots {_nT}_r$
tends to infinity and
$$
\frac{\phi({_nT})}{{_nT}_1\cdots {_nT}_r}\ \ge\ \frac{\phi(S)}{S_1\cdots S_r}
$$
for every $S\le {_nT}$. In particular, one can choose $S=({_nT})^a$. Then
$$
\frac{\phi({_nT})}{{_nT}_1\cdots {_nT}_r}\ \le\
\mu^{-r(j+1)}\frac{\phi_j({_nT})}{{_nT}_1\cdots {_nT}_r (\log {_nT}_1)^{j+1}\cdots (\log
{_nT}_r)^{j+1}}
$$ $$
\hspace{30pt} + \frac{\phi({_nT})}{{_nT}_1\cdots {_nT}_r}\frac 1{({_nT}_1\cdots
{_nT}_r)^{1-\mu}},
$$
so that
$$
\frac{\phi({_nT})}{{_nT}_1\cdots {_nT}_r}\ \le\ \frac{ \mu^{-r(j+1)}
\frac{\phi_j({_nT})}{{_nT}_1\cdots {_nT}_r (\log {_nT}_1)^{j+1}\cdots (\log {_nT}_r)^{j+1}}}
{1-\frac 1{({_nT}_1\cdots {_nT}_r)^{1-\mu}}}.
$$
Since the right hand side converges we reach a contradiction. This implies that
$$
L\df\limsup_{T\ra\infty} \frac{\phi(T)}{T_1\cdots T_r}\=\limsup_{T\ra\infty}
\frac{\phi(T^\mu)}{T_1^\mu\cdots T_r^\mu}
$$
is finite. We get
\begin{eqnarray*}
L &\le& \mu^{-r(j+1)} + L\ \limsup_{T\ra\infty} \frac 1{(T_1\cdots T_r)^{1-\mu}}\\
&=& \mu^{-r(j+1)}.
\end{eqnarray*}
Since $\mu$ is arbitrary we get $L\le 1$. The lemma follows.
\qed

We now finish the proof of the theorem. We have that
$\frac{\phi(T_1,\dots,T_r)}{T_1\cdots T_r}$ tends to $1$ as $T_k\ra\infty$. 
Since $a\in A^-$ has only eigenvalues of absolute value $<1$ on $\n$ it follows that
$0<\det(1-am\mid\n)<1$ for every $am\in A^-M$. Let $0<\eps<1$. Let $\psi_\eps(T)$ and
$\phi_\eps(T)$ denote the same sums as above extended only over those classes $[\ga]$ with
$1-\eps<\det(1-a_\ga m_\ga \mid \n)<1$

\begin{lemma}\label{compare}
As $T_1,\dots,T_r\ra\infty$ we have
$$
\frac{\phi(T_1,\dots,T_r)-\phi_\eps(T_1,\dots,T_r)}{T_1\cdots T_r}\ \ra\ 0,
$$
and
$$
\frac{\psi(T_1,\dots,T_r)-\psi_\eps(T_1,\dots,T_r)}{T_1\cdots T_r}\ \ra\ 0.
$$
\end{lemma}

\prf
Note that $\det(1-a_\ga m_\ga\mid\n)$ can only be $\le 1-\eps$ if $a_\ga$ is close to a wall in
$A^-$. In other words, $a_\ga$ has to lie in one of the segments
$$
S_\al^c\=\{ a\in A^-\mid a^{-\al}<c(\eps,\al)\}
$$
 for some $\al\in\Delta$ and some $c(\eps,\al)>1$. With $c=\max(c(\eps,\al),\al\in\Delta)$ it
follows that
$$
\phi(T)-\phi_\eps(T)\ \le\ \sum_{j=1,\dots,r}\phi(T_1,\dots, c,\dots ,T_r).
$$
The $k$th summand on the right hand side does not depend upon $T_k$. This implies the
first assertion. For the second note that
$$
\phi(T)-\phi_\eps(T)\ge \rez{1-\eps}(\psi(T)-\psi_\eps(T)).
$$
The lemma follows.
\qed

It follows that        
$$
\frac{\phi_\eps(T_1,\dots,T_r)}{T_1\cdots T_r}\ \ra\ 1\ \ \ {\rm as}\ T_j\ra\infty.
$$
For each $[\ga]$ in this sum we have
$$
1-\eps\ <\ \det(1-a_\ga m_\ga \mid \n)\ <\ 1,
$$
and so
$$
\frac{1-\eps}{\det(1-a_\ga m_\ga\mid\n)}\ <\ 1\ <\ \frac{1}{\det(1-a_\ga m_\ga\mid\n)}.
$$
Summing up we get
$$
\frac{\phi_\eps(T_1,\dots,T_r)}{T_1\cdots T_r}(1-\eps)\ <\ \frac{\psi_\eps(T_1,\dots,T_r)}{T_1\cdots
T_r}\ <\ \frac{\phi_\eps(T_1,\dots,T_r)}{T_1\cdots T_r}.
$$
Since $\frac{\phi_\eps(T_1,\dots,T_r)}{T_1\cdots T_r}$ tends to $1$ it follows that
$$
1-\eps\ \le\ \liminf_{T_j\ra\infty} \frac{\psi_\eps(T_1,\dots T_r)}{T_1\cdots T_r}\ \le\
\limsup_{T_j\ra\infty} \frac{\psi_\eps(T_1,\dots T_r)}{T_1\cdots T_r}\ \le\ 1.
$$
The second part of Lemma \ref{compare} gives
$$
1-\eps\ \le\ \liminf_{T_j\ra\infty} \frac{\psi(T_1,\dots T_r)}{T_1\cdots T_r}\ \le\
\limsup_{T_j\ra\infty} \frac{\psi(T_1,\dots T_r)}{T_1\cdots T_r}\ \le\ 1.
$$
Since $\eps$ is arbitrary the prime geodesic theorem follows.
\qed

We will finish this section with a conjecture. For each $[\ga]\in\CE(\Ga)$ pick a closed
geodesic $g_\ga$ which is closed by $\ga$. Note that
$\la_\ga$ equals the volume of the unique maximal flat that contains the geodesic $g_\ga$. The
contribution of all geodesics that lie in a given flat $F$ grows like $\frac{(\log T_1)\cdots(\log T_r)}{\vol(F)}$, as
$T_1,\dots,T_r\ra\infty$. This motivates the following conjecture.

\begin{conjecture}
Let
$$
\pi(T_1,\dots,T_r)\df \#\{[\ga]\in\CE(\Ga)\mid a_\ga^{\al_1}\le T_1,\dots, a_\ga^{\al_r}\le T_r\}.
$$
Then, as $T_1,\dots,T_r\ra\infty$,
$$
\pi(T_1,\dots,T_r)\ \sim\ \frac{T_1}{\log T_1}\cdots \frac{T_r}{\log T_r}.
$$
\end{conjecture}

It is not hard to show that the conjecture holds for products of rank one spaces.

\appendix
\section{Appendix: An application to class numbers}
In this section we give a new asymptotic formula for class numbers of orders in number fields. It
is quite different from known results like Siegel's Theorem (\cite{ayoub}, Thm 6.2). The
asymptotic is in several variables and thus contains more information than a single variable one.
In a sense it states that the units of the orders are equally distributed in different directions if only
one averages over sufficiently many orders.

Let
$d$ be a prime number
$\ge 3$. 
Let
$S$ be a finite set of primes with
$|S|\ge 2$. Let
$C(S)$ be the set of all totally real number fields $F$  with the property
$p\in S\
\ \Rightarrow\ \ p$ is non-decomposed in $F$.

Let $O(S)$ denote the set of all orders $\CO$ in number fields $F\in
C(S)$ which are maximal at each $p\in S$. For such an order $\CO$ let $h(\CO)$ be
its class number,
$R(\CO)$ its regulator and $\la_S(\CO)=\prod_{p\in S} f_p$, where $f_p$ is the inertia degree of
$p$ in
$F=\CO\otimes\Q$. Then $f_p\in \{ 1,d\}$ for every $p\in S$.

For $\la\in\CO^\times$ let $\rho_1,\dots,\rho_{d}$ denote the real embeddings of $F$ ordered
in a way that $|\rho_k(\la)|\ge |\rho_{k+1}(\la)|$ holds for $k=1,\dots,d-1$.

For $k=1,\dots, d-1$ let
$$
\al_k(\la)\df k(d-k)\,\log\left(
\frac{|\rho_{k}(\la)|}{|\rho_{k+1}(\la)|}\right).
$$

Let
$$
c=(\sqrt 2)^{1-d}  \left(\prod_{k=1}^{d-1}2k(d-k)\right).
$$
 So $c>0$ and it comes about as
correctional factor between the Haar measure normalization used in the Prime Geodesic
Theorem and the normalization used in the definition of the regulator.

\begin{theorem}\label{A.1}
For $T_1,\dots,T_r>0$ set
$$
\vartheta_S(T)\df \sum_{\stackrel{\la\in\CO^\times/\pm
1,\ \CO\in O(S)}{\stackrel{0<\al_k(\la)\le T_k}{k=1,\dots,d-1}}}  R(\CO)\,
h(\CO)\,
\la_S(\CO).
$$
Then we have, as $T_1,\dots,T_{d-1}\ra \infty$,
$$
\vartheta(T_1,\dots,T_{d-1})\ \sim\ \frac{c}{\sqrt{d}}\,T_1\cdots T_{d-1}.
$$
\end{theorem}

\prf
For given
$S$ there is a division algebra
$M$ over
$\Q$ of degree
$d$ which splits exactly outside $S$. Fix a maximal order $M(\Z)$ in $M$ and for
any ring $R$ define
$M(R)\df M(\Z)\otimes R$. Let $\det :M(R)\ra R$ denote the reduced norm then
$$
\CG(R)\df\{  x\in M(R)\mid \det(x)=1\}
$$
defines a group scheme over $\Z$ with $\CG(\R)\cong\SL_d(\R)=G$. Then $\Ga= \CG(\Z)$ is a
cocompact discrete torsion-free subgroup of $G$ (see \cite{class}). As can be seen in
\cite{class}, Theorem \ref{A.1} can be deduced from the prime geodesic theorem.
\qed

{\small University of Exeter, Mathematics, Exeter EX4
4QE, England\\ a.h.j.deitmar@ex.ac.uk}

\end{document}